\documentclass[12pt]{amsart}
\usepackage{amssymb}
\usepackage[all]{xy}

\newcommand{\issuenumber}{15}
\newcommand{\issuemonth}{December}
\newcommand{\issueyear}{2005}

\newtheorem{thm}{Theorem}[section]
\newtheorem{prob}[thm]{Problem}

\newtheorem{issue}{Issue}

\theoremstyle{definition}

\theoremstyle{remark}

\newcommand{\ed}{\end{thebibliography}\general\end{document}}

\newcommand{\fd}{\mathfrak{d}}
\newcommand{\fp}{\mathfrak{p}}

\newcommand{\NON}{{\mathsf   {NON}}}
\newcommand{\COF}{{\mathsf   {COF}}}

\newcommand{\cA}{\mathcal{A}}
\newcommand{\M}{\mathcal{M}}

\newcommand{\op}{\operatorname}
\newcommand{\cov}{\mathsf{cov}}

\newcommand{\cf}{\mathsf{cf}}

\newcommand{\R}{\mathbb{R}}

\newcommand{\fo}{\mathfrak{od}}

\renewcommand{\b}{\mathfrak{b}}

\renewcommand{\split}{\mathsf{Split}}
\newcommand{\bq}{\begin{quote}}
\newcommand{\eq}{\end{quote}}
\renewcommand{\O}{\mathcal{O}}
\newcommand{\B}{\mathcal{B}}
\newcommand{\BG}{\B_\Gamma}

\newcommand{\sone}{\mathsf{S}_1}    \newcommand{\sfin}{\mathsf{S}_{fin}}
\newcommand{\Sc}{\mathsf{S}_c}
\newcommand{\ufin}{\mathsf{U}_{fin}}

\newcommand{\nin}{\not\in}

\newcommand{\naturals}{{\mathbb N}}
\newcommand{\N}{\naturals}

\newcommand{\by}[2]{\par\hfill\emph{#1}, #2}
\newcommand{\nby}[1]{\par\hfill\emph{#1}}
\newcommand{\Tau}{\mathrm{T}}
\newcommand{\CE}{\textsc{CE}}

\newcommand{\be}{\begin{enumerate}}
\newcommand{\ee}{\end{enumerate}}
\newcommand{\bi}{\begin{itemize}}
\newcommand{\ei}{\end{itemize}}

\newcommand{\general}{\small\vfill\par\noindent\hrulefill\par
\noindent\textbf{Previous issues.} The first issues of this
bulletin, which contain general information (first issue), basic
definitions, research announcements, and open problems (all
issues) are available online, on \arx{math.GN/$x$}, where $x$ is
\texttt{0301011}, \texttt{0302062}, \texttt{0303057},
\texttt{0304087}, \texttt{0305367}, \texttt{0312140},
\texttt{0401155}, \texttt{0403369}, \texttt{0406411},
\texttt{0409072}, \texttt{0412305}, \texttt{0503631},
\texttt{0508563}, and \texttt{0509432},
respectively, for issues number $1$ to $14$.\\[0.1cm]
\textbf{Contributions.}
Please submit your contributions (announcements, discussions, and open problems)
by e-mailing us. It is preferred to write them
in \LaTeX{}.
The authors are urged to use as standard notation as possible, or otherwise give
the definitions or a reference to where the notation is explained.
Contributions to this bulletin would not require any transfer of copyright,
and material presented here can be published elsewhere.\\[0.1cm]
\textbf{Subscription.}
To receive this bulletin (free) to your
e-mailbox, e-mail us:\\
{boaz.tsaban@weizmann.ac.il}
}

\newcommand{\nArxPaper}[5]{\subsection{#2}{#4}\par\hfill{\arx{#1}}\par\hfill\emph{#3}}

\newcommand{\arx}[1]{\texttt{http://arxiv.org/abs/#1}}
\newcommand{\url}[1]{\bq\texttt{#1}\eq}
\newcommand{\online}[1]{The paper is available online at \url{#1}}

\title[$\mathcal{SPM}$ Bulletin \textbf{\issuenumber} (\issuemonth{} \issueyear)]{%
$\mathcal{SPM}$ Bulletin\\[0.5cm]
Issue number \issuenumber: \issuemonth{} \issueyear{} \CE{}}

\begin{document}
\maketitle


\section{Editor's note}

This issue is published while the very the feeling of the nearby Second SPM Workshop
is in the air. More than 50 mathematicians are expected to participate in this
fascinating event. Regardless whether you could or could not make it to this meeting,
you are encouraged to attend the coming BEST meeting (details below), in which a significant
part deals with SPM and related areas.

\medskip

Contributions to the next issue are, as always, welcome.

\medskip

\by{Boaz Tsaban}{boaz.tsaban@weizmann.ac.il}

\hfill \texttt{http://www.cs.biu.ac.il/\~{}tsaban}

\section{On Selective screenability and examples of R.\ Pol}

Let $\mathcal{A}$ and $\mathcal{B}$ be families of collections of
subsets of the infinite set $S$. The symbol
$\Sc(\mathcal{A},\mathcal{B})$ denotes the statement: For each
sequence $(O_n:n\in\omega)$ of elements of $\mathcal{A}$ there  is
a sequence $(T_n:n<\omega)$ where for each $n$ the family $T_n$ is
pairwise disjoint, and each element of $T_n$ is a subset of some
element of $O_n$, and $\cup_{n<\omega}T_n$ is an element of
$\mathcal{B}$. The author introduced this selection principle.

For $X$ a topological space the symbol $\mathcal{O}$ denotes the
collection of all open covers of $X$. The property
$\Sc(\mathcal{O},\mathcal{O})$ was introduced by Addis and Gresham
(1978). It is a selective version of Bing's property of
screenability (1951). Spaces with property
$\Sc(\mathcal{O},\mathcal{O})$ are traditionally called C-spaces.
According to a theorem of Hattori-Yamada, and independently Rohm
\cite{Rohm}, if $X$ is $\sigma$-compact and both $X$ and $Y$ have
$\Sc(\mathcal{O},\mathcal{O})$, then their product has
$\Sc(\mathcal{O},\mathcal{O})$. Consequently, if $X$ is
$\sigma$-compact and has property $\Sc(\mathcal{O},\mathcal{O})$,
then it has $\Sc(\mathcal{O},\mathcal{O})$ in all finite powers.

$\sigma$-compactness is stronger than the Hurewicz property, which
in turn is stronger than the Menger property \cite{SFH}. For
definitions of these properties see \cite{coc2}. R. Pol \cite{Pol}
showed that in the above product theorem the $\sigma$-compactness
of $X$ cannot be weakened to the Menger property, by proving in
Theorem 1.5 of \cite{Pol}: Assuming the Continuum Hypothesis, for
each natural number $n\ge 1$ there exists a separable metrizable
space $X$ such that: \be
\item $X^{n+1}$ has the Menger Property, and\\
\item $X^n$ is a C-space, but
\item $X^{n+1}$ is not a C-space.
\ee

We prove that the Menger property cannot be strengthened to the Hurewicz property in such examples:
\begin{thm}
Let $X$ be a separable metrizable space with property
$\Sc(\mathcal{O},\mathcal{O})$. For each $n$, if $X^n$ has the
Hurewicz property, then $X^n$ has property
$\Sc(\mathcal{O},\mathcal{O})$.
\end{thm}

Thus, if $X$ has property $\Sc(\mathcal{O},\mathcal{O})$, and if
it has the Hurewicz property in all finite powers, then it has the
property $\Sc(\mathcal{O},\mathcal{O})$ in all finite powers.

\nby{Liljana Babinkostova}

\section{Workshops and conferences}

\subsection{The Oxford Conference on Topology and Computer Science in Honour of
Peter Collins and Mike Reed}
7--10 August 2006, Oxford, UK.
This conference will mark the retirements of Peter Collins and Mike Reed from the
Oxford Mathematical Institute and the Oxford Computing Laboratory.
The conference will take the week before the Prague Symposium and will
include special sessions on
Set-theoretic Topology, Theoretical Computer Science and Continuum Theory
and Dynamics.

\subsection{Boise Extravaganza In Set Theory (BEST2006)}
March 31--April 2, 2006,
Boise State University
Boise, ID, USA.

We are pleased to announce our fifteenth annual BEST conference.
There will be four talks by invited speakers:
Natasha Dobrinen (Kurt Godel Research Center for Mathematical Logic),
Michael Hrusak (UNAM),
Istvan Juhasz (Alfred Renyi Institute of Mathematics, Budapest),
Boban Velickovic (Equipe de Logique Mathematique, Universite de Paris 7).

The current list of participants consists of
Liljana Babinkostova (Boise State University),
Randall Holmes (Boise State University),
Richard Laver (University of Colorado, Boulder),
Justin Moore (Boise State University),
Marion Scheepers (Boise State University), and
Boaz Tsaban (Weizmann Institute of Science).

If you wish to participate but are not listed above or if you are not on
the mailing list for the conference and wish to be, please contact the
organizers.

The conference webpage at
\url{http://diamond.boisestate.edu/~best/best15/best15.html}
contains the most current information including lodging, abstract
submission, travel, schedule, etc.

The talks will take up a full day
on Friday and Saturday and a half day (the morning) on Sunday. The exact
details of the schedule are not yet available. If you wish to participate
but will not attend the entire conference, please let us know in advance.

Financial Support: Limited financial support is available for
participants - particularly graduate students and young researchers -
who do not already have sufficient support to attend the conference. In
order to apply, e-mail one of the organizers.

Abstracts (generously provided by Atlas Conferences Inc.)
\url{http://atlas-conferences.com/cgi-bin/abstract/casb-01}

Organized by Liljana Babinkostova, Stefan Geschke, Justin Moore, and
Marion Scheepers.
You can reach us by writing to
\texttt{best@math.boisestate.edu}

The conference is supported by a grant from the National Science
Foundation (NSF grant DMS 0503444), whose assistance is gratefully
acknowledged.

\section{Research announcements}

\nArxPaper{math.GN/0509402}
{The isometry group of the Urysohn space as a Levy group}
{Vladimir Pestov}
{\newcommand{\Ur}{\mathbb{U}}\newcommand{\Iso}{\op{Iso}}
We prove that the isometry group $\Iso(\Ur)$ of the universal Urysohn metric
space $\Ur$ equipped with the natural Polish topology is a L\'evy group in the
sense of Gromov and Milman, that is, admits an approximating chain of compact
(in fact, finite) subgroups, exhibiting the phenomenon of concentration of
measure. This strengthens an earlier result by Vershik stating that $\Iso(\Ur)$
has a dense locally finite subgroup. We propose a reformulation of Connes'
Embedding Conjecture as an approximation-type statement about the unitary group
$U(\ell^2)$, and show that in this form the conjecture makes sense also for
$\Iso(\Ur)$.}

\nArxPaper{math.LO/0509392}
{Chasing Silver}
{Andrzej Roslanowski and Juris Steprans}
{We show that limits of CS iterations of $n$-Silver forcing notion have the
$n$-localization property.}

\nArxPaper{math.LO/0509406}
{Disjoint Non-Free Subgoups of Abelian Groups}
{Andreas Blass and Saharon Shelah}
{Let $G$ be an abelian group and let $\lambda$ be the smallest rank of any group
whose direct sum with a free group is isomorphic to $G$. If $\lambda$ is
uncountable, then $G$ has $\lambda$ pairwise disjoint, non-free subgroups. There is
an example where $\lambda$ is countably infinite and $G$ does not have even two
disjoint, non-free subgroups.}

\nArxPaper{math.GT/0509395}
{Characterizing metric spaces whose hyperspaces are absolute neighborhood retracts}
{T.\ Banakh and R.\ Voytsitskyy}
{We characterize metric spaces $X$ whose hyperspaces $2^X$ or $Bd(X)$ of
non-empty closed (bounded) subsets, endowed with the Hausdorff metric, are
absolute [neighborhood] retracts.}

\subsection{A Vitali set can be homeomorphic to its complement}
We prove that:
\begin{enumerate}
\item There exists a special automorphism $F$ of the Cantor
space such that every noncancellable composition of finite
powers of $F$ and translations of rational numbers has no
fixed point.
\item For this automorphism there exists both
Vitali and Bernstein subset of the Cantor space
such that the image by $F$ of this set is equal to its
complement.
\item There exists a Bernstein and Vitali set
such that there is no Borel isomorphism between
this set and its complement.
\end{enumerate}
\par\hfill{\texttt{http://delta.univ.gda.pl/\~{}andrzej/vitali.zip}}

\nby{Andrzej Nowik}

\nArxPaper{math.GN/0509607}
{$o$-Boundedness of free objects over a Tychonoff space}
{Lyubomyr Zdomskyy}
{We characterize  various sorts of boundedness of the free  (abelian)
topological group $F(X)$ ($A(X)$) as well as the free locally-convex linear topological space
$L(X)$ in terms of properties of a Tychonoff space $X$.
These properties appear
to be close to so-called selection principles, which permits us to
show, that (it is consistent with ZFC that) the property of Hurewicz (Menger) is $l$-invariant.
This gives a method of construction of
$OF$-undetermined topological groups with strong combinatorial properties.
}

\subsection{On the consistency strength of the Milner-Sauer conjecture}
Motivated by a conjecture of Milner and Sauer
regarding partial orders of singular cofinality,
we came to prove a theorem about topological spaces of singular density.

We prove that the existence of a topological space $\langle X,O\rangle$
satisfying $d(X)=w(X)=\lambda>\cf(\lambda)\ge\hat{h}(X)$ implies
$\cf([\lambda]^{<\cf(\lambda)},\subseteq)>\lambda$.
One consequence is that the existence of a hereditarily Lindel\"of
topological space of density and weight $\aleph_{\omega_1}$ implies the
existence of a measurable cardinal in an inner model of ZFC.

On the way, we also notice that for any topological space
$\langle X,O \rangle$, if $d(X)$ is a singular cardinal,
then $|O|>d(X)$.

\hfill\texttt{http://dx.doi.org/10.1016/j.apal.2005.09.012}

\nby{Assaf Rinot}

\nArxPaper{math.GN/0510118}
{Reconstruction of manifolds and subsets of normed spaces from subgroups of their homeomorphism groups}
{Matatyahu Rubin and Yosef Yomdin}
{We prove various reconstruction theorems about open subsets of normed spaces.
E.g. if the uniformly continuous homeomorphism groups of two such sets are
isomorphic, then this isomorphism is induced by a uniformly continuous
homeomorphism between these open sets.}

\nArxPaper{math.GN/0510120}
{Reconstruction theorem for homeomorphism groups without small sets and non-shrinking functions of a normed space}
{Vladimir P.\ Fonf and Matatyahu Rubin}
{We prove a reconstruction theorem for homeomorphism groups of open sets in
metrizable locally convex topological vector spaces. We show that certain small
subgroups of the full homeomorphism group obey the conditions of the above
theorem.}

\nArxPaper{math.LO/0510122}
{Locally Moving Groups and the Reconstruction Problem for Chains and Circles}
{Stephen McCleary and Matatyahu Rubin}
{We prove that complete Boolean algebras can be reconstructed from any locally
moving subgroup of their full automorphism group. We use this theorem in order
to prove that linear orders and circles can be reconstructed from small
subgroups of their full automorphism groups.}

\nArxPaper{math.CO/0510254}
{Divisibility of countable metric spaces}
{Christian Delhomme, Claude Laflamme, Maurice Pouzet, and Norbert Sauer}
{Prompted by a recent question of G.\ Hjorth as to whether
a bounded Urysohn space is indivisible, that is to say has the
property that any partition into finitely many pieces has one piece
which contains an isometric copy of the space, we answer this
question and more generally investigate partitions of countable metric
spaces.

We show that an indivisible metric space must be totally Cantor
disconnected, which implies in particular that every Urysohn space $\mathbb{U}_{V}$ with $V$
bounded or not but dense in some initial segment of $\R_+$, is divisible.
On the other hand we also show that one
can remove ``large'' pieces from a bounded Urysohn space with the
remainder still inducing a copy of this space, providing a certain
``measure'' of the indivisibility.  Associated with every totally
Cantor disconnected space is an ultrametric space, and we go on to
characterize the countable ultrametric spaces which are homogeneous
and indivisible.
}

\nArxPaper{math.FA/0510407}
{Pre-compact families of finite sets of integers and weakly null sequences in Banach spaces}
{Jordi Lopez Abad and Stevo Todorcevic}
{We provide a somewhat general framework for studying weakly null sequences in
Banach spaces using Ramsey theory of families of finite subsets of integers.}

\nArxPaper{math.GN/0510484}
{A semifilter approach to selection principles II: $\tau^*$-covers}
{Lyubomyr Zdomskyy}
{In this paper we settle all questions whether the
properties $P$ and $Q$ provably coincide, where $P$ and $Q$ run over selection
principles of the type $\ufin(\O,\cA)$.}

\nArxPaper{math.LO/0510477}
{Parametrizing the abstract Ellentuck theorem}
{Jose G.\ Mijares}
{We give a parametrization with perfect sets of the abstract Ellentuck
theorem. The main tool for achieving this goal is a sort of parametrization of
an abstract version of the Nash-Williams theorem. As corollaries, we obtain
some known classical results like the parametrized version of the Galvin-Prikry
theorem due to Miller and Todorcevic, and the parametrized version of
Ellentuck's theorem due to Pawlikowski. Also, we obtain parametized vesions of
nonclassical results such as Milliken's theorem.
}

\nArxPaper{math.LO/0510517}
{A notion of selective ultrafilter corresponding to topological Ramsey spaces}
{Jose G.\ Mijares}
{We introduce the relation of \textit{almost-reduction} in an arbitrary topological
Ramsey space $\mathcal{R}$ as a generalization of the relation of
almost-inclusion on $\N^{[\infty]}$. This leads us to a type of
ultrafilter $\mathcal{U}$ on the set of first approximations of
the elements of $\mathcal{R}$ which corresponds to the well-known
notion of \textit{selective ultrafilter} on $\N$. The relationship
turns out to be rather exact in the sense that it permits us to
lift several well-known facts about selective ultrafilters on
$\mathbb{N}$ and the Ellentuck space  $\mathbb{N}^{[\infty]}$ to
the ultrafilter $\mathcal{U}$ and the Ramsey space $\mathcal{R}$.
For example, we prove that the Open Coloring Axiom holds on
$L(\mathbb{R})[\mathcal{U}]$, extending therefore the result from
[Prisco and Todorcevic, \emph{Perfect-set properties in $L(\mathbb{R})[\mathcal{U}]$},
Adv.\ Math. \textbf{139} (1998), 240-259],
which gives the same conclusion for the Ramsey space
$\N^{[\infty]}$.}

\nArxPaper{math.GN/0511567}
{Compact spaces generated by retractions}
{Wieslaw Kubis}
{We study compact spaces which are obtained from metric compacta by iterating
the operation of inverse limit of continuous sequences of retractions. We
denote this class by R. Allowing continuous images in the definition of class
R, one obtains a strictly larger class, which we denote by RC. We show that
every space in class RC is either Corson compact or else contains a copy of the
ordinal segment $[0,\aleph_1]$. This improves a result of Kalenda, where the
same was proved for the class of continuous images of Valdivia compacta. We
prove that spaces in class R do not contain cutting P-points (see the
definition below), which provides a tool for finding spaces in RC minus R.
Finally, we study linearly ordered spaces in class RC. We prove that scattered
linearly ordered compacta belong to RC and we characterize those ones which
belong to R. We show that there are only 5 types (up to order isomorphism) of
connected linearly ordered spaces in class R and all of them are Valdivia
compact. Finally, we find a universal pre-image for the class of all linearly
ordered Valdivia compacta.}

\nArxPaper{math.GN/0511437}
{Gromov-Hausdorff ultrametric}
{Ihor Zarichnyi}
{We show that there exists a natural counterpart of the Gromov-Hausdorff
metric in the class of ultrametric spaces. It is proved, in particular, that
the space of all ultrametric spaces whose metric take values in a fixed
countable set is homeomorphic to the space of irrationals.
}

\nArxPaper{math.FA/0511456}
{Computing the complexity of the relation of isometry between separable Banach spaces}
{Julien Melleray}
{We compute here the Borel complexity of the relation of isometry between
separable Banach spaces, using results of Gao, Kechris and Weaver.}

\nArxPaper{math.GN/0511606}
{On some classes of Lindel\"of Sigma-spaces}
{Wieslaw Kubis, Oleg Okunev, and Paul J.\ Szeptycki}
{We consider special subclasses of the class of Lindel\"of Sigma-spaces
obtained by imposing restrictions on the weight of the elements of compact
covers that admit countable networks: A space $X$ is in the class
$L\Sigma(\leq\kappa)$ if it admits a cover by compact subspaces of weight
$\kappa$ and a countable network for the cover. We restrict our attention to
$\kappa\leq\omega$. In the case $\kappa=\omega$, the class includes the
class
of metrizably fibered spaces considered by Tkachuk, and the $P$-approximable
spaces considered by Tkacenko. The case $\kappa=1$ corresponds to the spaces
of
countable network weight, but even the case $\kappa=2$ gives rise to a
nontrivial class of spaces. The relation of known classes of compact spaces
to
these classes is considered. It is shown that not every Corson compact of
weight $\aleph_1$ is in the class $L\Sigma(\leq \omega)$, answering a
question
of Tkachuk. As well, we study whether certain compact spaces in
$L\Sigma(\leq\omega)$ have dense metrizable subspaces, partially answering a
question of Tkacenko. Other interesting results and examples are obtained,
and
we conclude the paper with a number of open questions.
}

\nArxPaper{math.LO/0512217}
{On the depth of Boolean algebras}
{Saharon Shelah}
{We show that the $\op{Depth}^+$ of an ultraproduct of Boolean algebras, can not
jump over the $\op{Depth}^+$ of every component by more than one cardinality. We can
have, consequently, similar results for the $\op{Depth}$ invariant.
}

\section{Problem of the Issue}

In the paper
\url{http://arxiv.org/abs/math.GN/0510162}
we formulate three problems concerning topological properties
of sets generating Borel non-$\sigma$-compact groups. In case of the concrete
$F_{\sigma\delta}$-subgroup of the Cantor group this gives an equivalent
reformulation of the Scheepers diagram problem.
The problems are related to the following open problem.

\begin{prob}
Can a Borel non-$\sigma$-compact group be generated by a Hurewicz subspace?
\end{prob}

\nby{Lyubomyr Zdomskyy}

\section{Problems from earlier issues}

\begin{issue}
Is $\binom{\Omega}{\Gamma}=\binom{\Omega}{\Tau}$?
\end{issue}

\begin{issue}
Is $\ufin(\Gamma,\Omega)=\sfin(\Gamma,\Omega)$?
And if not, does $\ufin(\Gamma,\Gamma)$ imply
$\sfin(\Gamma,\Omega)$?
\end{issue}

\stepcounter{issue}

\begin{issue}
Does $\sone(\Omega,\Tau)$ imply $\ufin(\Gamma,\Gamma)$?
\end{issue}

\begin{issue}
Is $\fp=\fp^*$? (See the definition of $\fp^*$ in that issue.)
\end{issue}

\begin{issue}
Does there exist (in ZFC) an uncountable set satisfying $\sone(\BG,\B)$?
\end{issue}

\stepcounter{issue}

\begin{issue}
Does $X \nin \NON(\M)$ and $Y\nin\mathsf{D}$ imply that
$X\cup Y\nin \COF(\M)$?
\end{issue}

\begin{issue}
Assume CH. Is $\split(\Lambda,\Lambda)$ preserved under finite unions?
\end{issue}

\begin{issue}
Is $\cov(\M)=\fo$? (See the definition of $\fo$ in that issue.)
\end{issue}

\begin{issue}
Does $\sone(\Gamma,\Gamma)$ always contain an element of cardinality $\b$?
\end{issue}

\begin{issue}
Could there be a Baire metric space $M$ of weight $\aleph_1$ and a partition
$\mathcal{U}$ of $M$ into $\aleph_1$ meager sets where for each ${\mathcal U}'\subset\mathcal U$,
$\bigcup {\mathcal U}'$ has the Baire property in $M$?
\end{issue}

\stepcounter{issue} 

\begin{issue}
Does there exist (in ZFC) a set of reals $X$ of cardinality $\fd$ such that all
finite powers of $X$ have Menger's property $\ufin(\O,\O)$?
\end{issue}

\begin{issue}
Can a Borel non-$\sigma$-compact group be generated by a Hurewicz subspace?
\end{issue}

\begin{thebibliography}{00}

\bibitem{Sc} L.\ Babinkostova, \emph{Selective versions of screenability}, {\bf Filomat 17:2} (2003), 127 - 134

\bibitem{coc2} W.\ Just, A.W.\ Miller, M.\ Scheepers and P.J.\ Szeptycki, \emph{The combinatorics of open covers II}, {\bf Top. Appl.} 73 (1996), 241 -- 266.

\bibitem{Pol} R.\ Pol, \emph {Weak infinite-dimensionality in Cartesian products with the Menger property}, {\bf Topology and its Applications} 61 (1995), 85--94

\bibitem{Rohm} D.\ Rohm, \emph {Products of infinite-dimensional spaces}, {\bf Proceedings of the Amer. Math. Society} 108:4 (1990), 1019--1023

\bibitem{SFH}
B.\ Tsaban and L.\ Zdomsky,
\emph{Scales, fields, and a problem of Hurewicz}, submitted.\\
\arx{math.GN/0507043}

\ed